\documentclass[english,aps,prl,reprint,twocolumn]{revtex4-1}
\usepackage{prettyref}
\usepackage{amsmath}
\usepackage{amsthm}
\usepackage{amssymb}
\usepackage{graphicx}
\usepackage{xcolor}
\usepackage{tikz}
\usepackage{xcolor}
\usepackage{bbm}

\usepackage[normalem]{ulem}

\newcommand{\R}{\mathbb{R}}

\newcommand{\x}{\mathbf{x}}
\newcommand{\f}{\mathbf{f}}
\newcommand{\ov}{\overline}

\newcommand\solidrule[1][21pt]{\rule[0.5ex]{#1}{.4pt}}
\newcommand\dashedrule{\mbox{%
	\solidrule[5pt]\hspace{3pt}\solidrule[5pt]\hspace{3pt}\solidrule[5pt]}}
\newcommand\dottedrule{\mbox{%
	\solidrule[2pt]\hspace{2pt}\solidrule[2pt]\hspace{2pt}\solidrule[2pt]\hspace{2pt}\solidrule[2pt]\hspace{2pt}\solidrule[2pt]\hspace{2pt}\solidrule[2pt]}}

\begin{document}

\title{Optimal bounds and extremal trajectories for\\time averages in nonlinear dynamical systems}

\author{Ian Tobasco$^1$}
\author{David Goluskin$^{1,4}$}
\author{Charles R.~Doering$^{1,2,3}$}

\affiliation{$^1$Department of Mathematics, University of Michigan, Ann Arbor, MI 48109, USA}
\affiliation{$^2$Department of Physics, University of Michigan, Ann Arbor, MI 48109, USA}
\affiliation{$^3$Center for the Study of Complex Systems, University of Michigan, Ann Arbor, MI 48109, USA}
\affiliation{$^4$Department of Mathematics and Statistics, University of Victoria, Victoria, BC, V8P 5C2, Canada}

\date{\today}
\begin{abstract}
For any quantity of interest in a system governed by ordinary differential equations, it is natural to seek the largest (or smallest) long-time average among solution trajectories, as well as the extremal trajectories themselves. Upper bounds on time averages can be proved \emph{a priori} using auxiliary functions, the optimal choice of which is a convex optimization problem. We prove that the problems of finding maximal trajectories and minimal auxiliary functions are strongly dual. Thus, auxiliary functions provide arbitrarily sharp upper bounds on time averages. Moreover, any nearly minimal auxiliary function provides phase space volumes in which all nearly maximal trajectories are guaranteed to lie. For polynomial equations, auxiliary functions can be constructed by semidefinite programming, which we illustrate using the Lorenz system.
\\
\\
{\bf Keywords:} Nonlinear dynamical systems, Time averages, Ergodic theory, Semidefinite programming, Sum-of-squares polynomials, Lorenz equations
\end{abstract}

\maketitle

\noindent
{\bf 1. Introduction} \quad
For dynamical systems governed by ordinary differential equations (ODEs) whose solutions are complicated and perhaps chaotic, the primary interest is often in long-time averages of key quantities. Time averages can depend on initial conditions, so it is natural to seek the largest or smallest average among all trajectories, as well as the extremal trajectories that realize them. For various purposes including the control of chaos \cite{Yang2000}, it is valuable to know extremal trajectories regardless of their stability. In other situations one is interested only in stable trajectories, but determining extrema only among these can be prohibitively difficult. The next best option is to determine extrema among all trajectories.

One common way to seek extremal time averages is to construct a large number of candidate trajectories. However, for many nonlinear systems it is challenging both to compute trajectories and to determine that the extremal ones have not been overlooked. In this Letter we study an alternative approach that is broadly applicable and often more tractable: constructing sharp \emph{a priori} bounds on long-time averages. We focus on upper bounds; lower bounds are analogous.

The search for an upper bound on a long-time average can be posed as a convex optimization problem \cite{Chernyshenko2014}, as described in the next section. Its solution requires no knowledge of trajectories. What is optimized is an auxiliary function defined on phase space, similar to but distinct from Lyapunov functions in stability theory. We prove here that the best bound produced by solving this convex optimization problem coincides exactly with the extremal long-time average. That is, arbitrarily sharp bounds on time averages can be produced using increasingly optimal auxiliary functions. Moreover, nearly optimal auxiliary functions yield volumes in phase space where maximal and nearly maximal trajectories must reside. Whether such auxiliary functions can be computed in practice depends on the system being studied, but when the ODE and quantity of interest are polynomial, auxiliary functions can be constructed by solving semidefinite programs (SDPs) \cite{Chernyshenko2014, Fantuzzi2016, Goluskin2017}. The resulting bounds can be arbitrarily sharp. We illustrate these methods using the Lorenz system \cite{Lorenz1963}.

Consider a well-posed autonomous ODE on $\R^d$,
\begin{equation}
\tfrac{d}{dt}\x=\f(\x),
\label{eq: ode}
\end{equation}
whose solutions are continuously differentiable in their initial conditions. To guarantee this, we assume that $\f(\x)$ is continuously differentiable. Given a continuous quantity of interest $\Phi(\x)$, we define its \emph{long-time average} along a trajectory $\x(t)$ with initial condition $\x(0)=\x_0$~by
\begin{equation}
\overline{\Phi}(\x_0)=\limsup_{T\to\infty}\,\frac{1}{T}\int_{0}^{T}\Phi(\mathbf{x}(t))\,dt.
\label{eq: Phi}
\end{equation}
Time averages could be defined using $\liminf$ instead; our results hold \emph{mutatis mutandis} \footnote{The $\limsup$ and $\liminf$ averages need not coincide on every trajectory, but their maxima over trajectories do.}.

Let $B\subset \R^d$ be a closed bounded region such that trajectories beginning in $B$ remain there. In a dissipative system $B$ could be an absorbing set; in a conservative system $B$ could be defined by constraints on invariants. We are interested in the maximal long-time average among all trajectories eventually remaining in $B$: 
\begin{equation}
\ov\Phi^*=\max_{\x_0\in B}\,\ov\Phi(\x_0).
\label{eq: Phi*}
\end{equation}
As shown below, there exist $\x_0$ attaining the maximum. The fundamental questions addressed here are: what is the value of $\ov\Phi^*$, and which trajectories attain it?

\smallskip
\noindent
{\bf 2. Bounds by convex optimization} \quad
Upper bounds on long-time averages can be deduced using the fact that time derivatives of bounded functions average to zero. Given any initial condition $\x_0$ in $B$ and any $V(\x)$ in the class $C^1(B)$ of continuously differentiable functions on $B$ \footnote{Here $C^1(B)$ denotes functions on $B$ admitting a continuously differentiable extension to a neighborhood of $B$.},
\begin{equation}
\ov{\tfrac{d}{dt}V}=\ov{\f\cdot\nabla V}=0.
\label{eq: identity}
\end{equation}
This generates an infinite family of functions with the same time average as $\Phi$ since for all such $V$
\begin{equation}
\ov\Phi=\ov{\Phi+\f\cdot\nabla V}.
\end{equation}
Bounding the righthand side pointwise gives
\begin{equation}
\ov\Phi(\x_0) \le \max_{\x\in B}\left\{ \Phi+\mathbf{f}\cdot\nabla V\right\}
\label{eq: weak 1}
\end{equation}
for all initial conditions $\x_0\in B$ and \emph{auxiliary functions} $V\in C^1(B)$.
Expression \eqref{eq: weak 1} is useful since no knowledge of trajectories is needed to evaluate the righthand side.

To obtain the optimal bound implied by \eqref{eq: weak 1}, we minimize the righthand side over $V$ and maximize the lefthand side over $\x_0$:
\begin{equation}
\max_{\x_0\in B}\, \ov\Phi
\le \inf_{V\in C^1(B)}
\max_{\x\in B}\left\{ \Phi+\mathbf{f}\cdot\nabla V\right\}.
\label{eq: weak 2}
\end{equation}
The minimization over auxiliary functions $V$ in \eqref{eq: weak 2} is convex, although minimizers need not exist. The main mathematical result of this Letter is that the lefthand and righthand optimizations are dual variational problems, and moreover that strong duality holds, meaning that \eqref{eq: weak 2} can be improved to an equality:
\begin{equation}
\max_{\x_0\in B}\, \ov\Phi
= \inf_{V\in C^1(B)} \max_{\x\in B}\left\{ \Phi+\mathbf{f}\cdot\nabla V\right\}.
\label{eq: strong}
\end{equation}
Thus, arbitrarily sharp bounds on the maximal time average $\ov\Phi^*$ can be obtained using increasingly optimal $V$.

The auxiliary function method is not the same as the various Lyapunov-type methods used to show stability or boundedness in ODE systems. However, in instances where $\Phi(\x)$ approaches infinity as $|\x|\to\infty$, auxiliary functions that imply finite upper bounds $\ov\Phi\le U$ also imply the existence of trapping sets by the following argument. Suppose  $V\in C^1(\mathbb{R}^d)$ is an auxiliary function for which the maximum of $\Phi+\mathbf{f}\cdot\nabla V$ over $\mathbb R^d$ is no larger than $U$. Then,
\begin{equation}
\tfrac{d}{dt}V=\f\cdot\nabla V \leq U - \Phi \to -\infty 
\label{eq: Lyap}
\end{equation}
as $|\x|\to\infty$. Expression \eqref{eq: Lyap} is a typical Lyapunov-type condition implying that all sufficiently large sublevel sets of $\Phi$ must be trapping sets.

The remainder of this Letter is organized as follows. The next section describes how nearly optimal $V$ can also be used to locate maximal and nearly maximal trajectories in phase space. The section after illustrates these ideas using the Lorenz system, for which we have constructed nearly optimal $V$ by solving SDPs. The final section proves the strong duality \eqref{eq: strong} and establishes the existence of maximal trajectories.

\smallskip
\noindent
{\bf 3. Near optimizers} \quad 
In light of the duality \eqref {eq: strong}, an initial condition $\x_0^*$ and auxiliary function $V^*$ are optimal if and only if they satisfy 
\begin{equation}
\ov\Phi(\x_0^*) = \max_{\x\in B}\left\{\Phi + \mathbf{f}\cdot \nabla V^* \right \}.
\end{equation}
Even if the infimum over $V$ in \eqref{eq: strong} is not attained, there exist nearly optimal pairs.
That is, for all $\epsilon>0$ there exist $(\x_0,V)$ for which \eqref{eq: weak 1} is within $\epsilon$ of an equality:
\begin{equation}
0\le \max_{\x\in B}\left\{ \Phi+\mathbf{f}\cdot\nabla V\right\} -  \ov\Phi(\x_0)  \le \epsilon.
\label{eq: near optimal def}
\end{equation}
In such cases, $\max_{\x\in B}\left\{ \Phi+\mathbf{f}\cdot\nabla V\right\}$ is within $\epsilon$ of being a sharp upper bound on $\ov\Phi^*$, while the trajectory starting at $\x_0$ achieves a time average $\ov\Phi$ within $\epsilon$ of $\ov\Phi^*$.

Nearly optimal $V$ can be used to locate all trajectories consistent with \eqref{eq: near optimal def}.
Moving the constant term inside the time average and subtracting the identity \eqref{eq: identity} gives
\begin{equation}
0\le \ov{\max_{\x\in B}\left\{ \Phi+\mathbf{f}\cdot\nabla V\right\}
	- (\Phi +\f\cdot\nabla V)} \le \epsilon
\label{eq: near optimal ineq}
\end{equation}
for such trajectories.
The integrand in \eqref{eq: near optimal ineq} is nonnegative, and the fraction of time it exceeds $\epsilon$ can be estimated. Consider the set where the integrand is no larger than~$M>\epsilon$,
\begin{equation}
\mathcal S_M = \Big\{ \x \in B: \max_{\x\in B}\left\{ \Phi+\mathbf{f}\cdot\nabla V\right\} - (\Phi+\mathbf{f}\cdot\nabla V)(\x)  \le M \Big\}.
\label{eq: S}
\end{equation}
Let $\mathcal F_M(T)$ denote the fraction of time $t\in[0,T]$ during which $\x(t)\in\mathcal S_M$.
For any trajectory obeying \eqref{eq: near optimal ineq}, this time fraction is bounded below as
\begin{equation}
\liminf_{T\to\infty}\, \mathcal F_M(T) \ge1-\epsilon/M.
\label{eq: Cheb est}
\end{equation}
This follows from an application of Markov's inequality: as the integrand in \eqref{eq: near optimal ineq} is nonnegative,
\begin{equation}
\epsilon \geq \ov{ M \mathbbm{1}_{\x\notin \mathcal S_M} } = M\left( 1 - \liminf_{T\to\infty}\, \mathcal F_M(T)\right ).
\end{equation}

In practice, it may not be known if there exist trajectories satisfying \eqref{eq: near optimal def} for a given $V$ and $\epsilon$. Still, the estimate \eqref{eq: Cheb est} says that any such trajectories would lie in $\mathcal S_M$ for a fraction of time no smaller than $1-\epsilon/M$. The conclusion is strongest when $\epsilon\ll M$, but if $M$ is too large the volume $\mathcal S_M$ is large and featureless, failing to distinguish nearly maximal trajectories. The result is most informative when $V$ is nearly optimal so that there exist trajectories where $\epsilon \ll M$ with $M$ not too large.

If a minimal $V^*$ exists, its set $\mathcal S_0$ is related to maximal trajectories. Any such trajectory achieves $\epsilon=0$ in \eqref{eq: near optimal ineq}. If it is a periodic orbit, for instance, it must lie in $\mathcal{S}_{0}$.
Thus $V^{*}$ is determined up to a constant on maximal orbits. More generally, $V^*$ must satisfy 
\begin{equation}
\Phi(\x) +\f(\x) \cdot \nabla V^*(\x) = \ov\Phi^*
\label{eq: optimality conds}
\end{equation}
for all $\x\in\mathcal{S}_0$.
It is tempting to conjecture that $\mathcal S_0$ coincides with maximal trajectories but, as described at the end of the next section, $\mathcal{S}_{0}$ can also contain points not on any maximal trajectory.

\smallskip
\noindent
{\bf 4. Nearly optimal bounds and orbits in the Lorenz system} \quad
When $\f(\x)$ and $\Phi(\x)$ are polynomials, $V(\x)$ can be optimized computationally within a chosen polynomial ansatz by solving an SDP \cite{Chernyshenko2014, Fantuzzi2016, Goluskin2017}. The bound $\ov\Phi^*\le U$ follows from \eqref{eq: weak 1} if $\Phi+\mathbf{f}\cdot\nabla V\le U$ for all $\x\in B$. A sufficient condition for this is that $U-\Phi-\mathbf{f}\cdot\nabla V$ is a sum of squares (SOS) of polynomials. The latter is equivalent to an SDP and is often computationally tractable \cite{Parrilo2003, Lasserre2010}.

It does not follow from the strong duality result \eqref{eq: strong} that bounds computed by SOS methods can be arbitrarily sharp. This is because requiring the polynomial $U-\Phi-\mathbf{f}\cdot\nabla V$ to be SOS is generally stronger than requiring it to be nonnegative \cite{Hilbert1888, Parrilo2003}. Nonetheless, in the few examples where time averages have been bounded using SOS methods \cite{Fantuzzi2016, Goluskin2017}, the bounds either are sharp or appear to become sharp as the polynomial degree of $V$ increases.

\begin{figure}
\centering
\begin{tikzpicture}
\node at (0,0) {\includegraphics[width=.98\columnwidth]{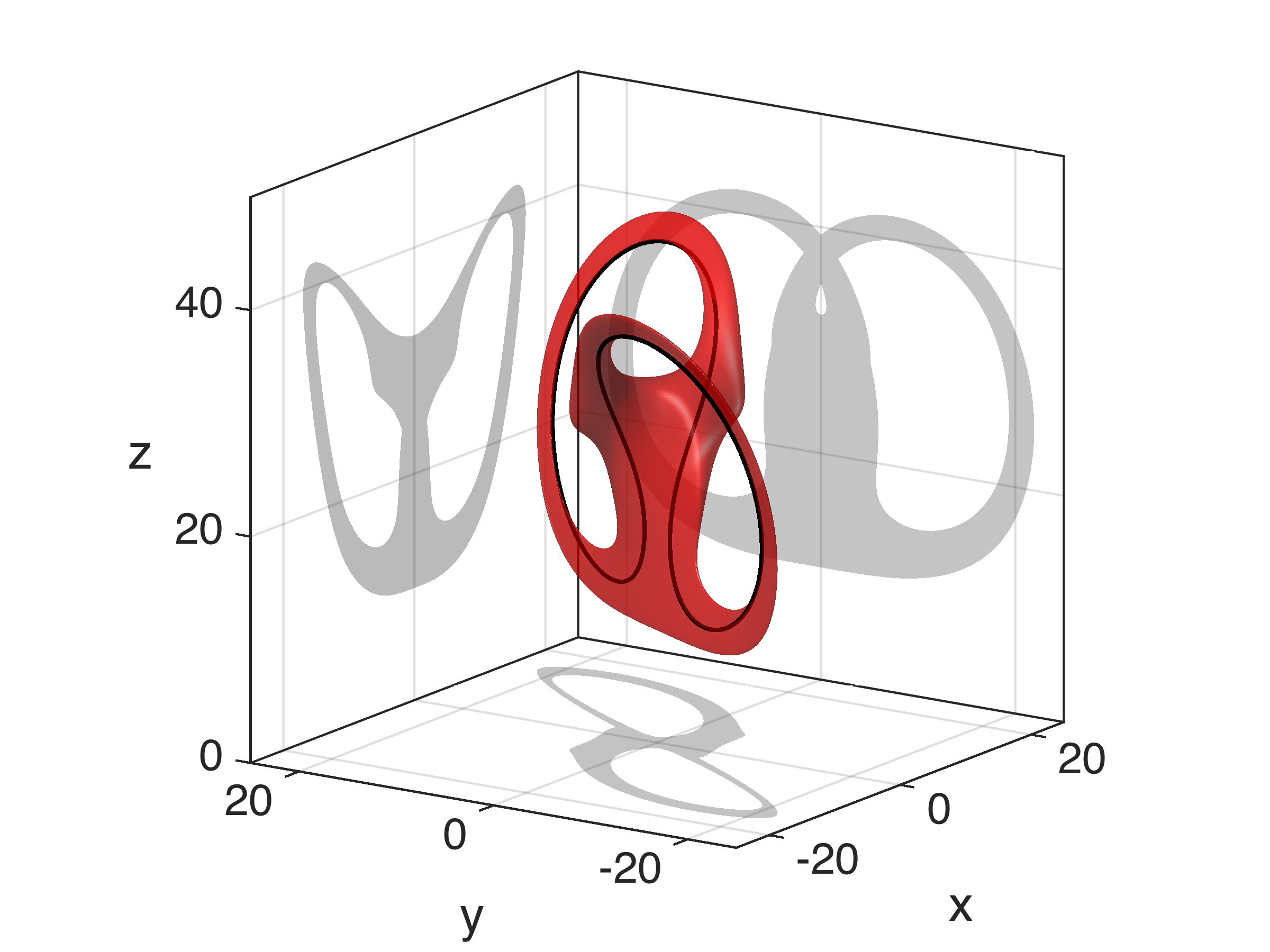}};
\node at (-3.33,2.6) {(a)};
\end{tikzpicture}\\
\begin{tikzpicture}
\node at (0,0) {\includegraphics[width=.98\columnwidth]{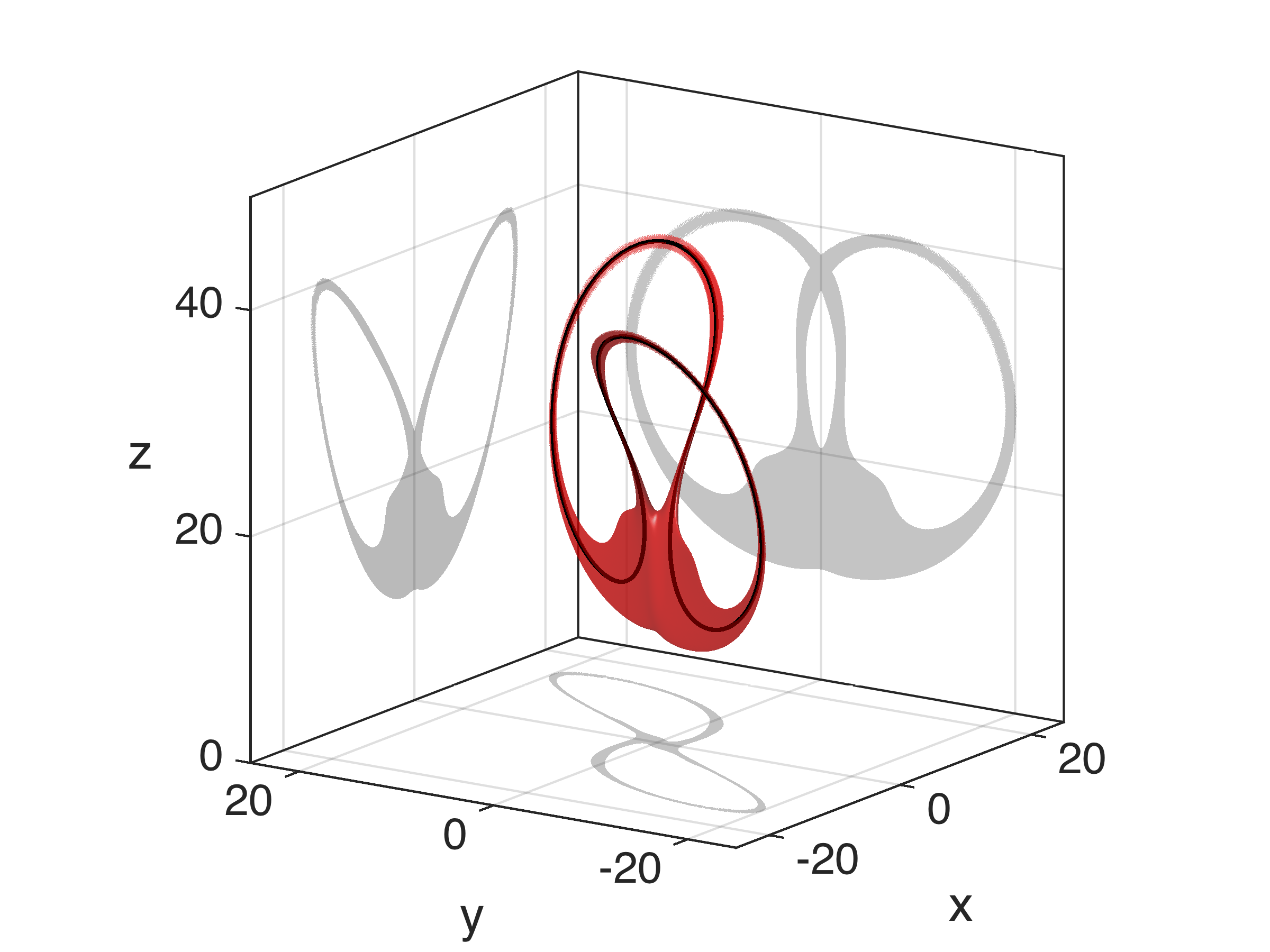}};
\node at (-3.33,2.6) {(b)};
\end{tikzpicture}
\caption{\label{fig: 3D}The volumes in phase space (a) $\mathcal S_{3000}$ for the optimal degree-6 polynomial $V$ and (b) $\mathcal S_{1000}$ for the optimal degree-10 polynomial $V$. Any trajectory maximizing $\ov{z^4}$ must spend at least $99.97\%$ of its time in $\mathcal S_{1000}$. The black curves show the shortest periodic orbit, which appears to maximize $\ov{z^4}$.}
\end{figure}

The remainder of this section presents the results of SOS bounding computations for the Lorenz system at the standard chaotic parameters $(\beta,\sigma,r)=(8/3,10,28)$. We obtain nearly sharp bounds on the maximal time average of $\Phi(x,y,z)=z^4$, as well as approximations to maximal trajectories. Because there exist compact absorbing balls \cite{Lorenz1979}, maximization over such $B$ in \eqref{eq: Phi*} is equivalent to maximization over $\R^d$. As reported in \cite{Goluskin2017}, searching among the periodic orbits computed by Viswanath \cite{Viswanath2003} suggests that the maximal average $\ov{z^4}^*$ is attained by the shortest periodic orbit---the black curves in Fig.~\ref{fig: 3D}. We have used SOS methods to construct nearly optimal $V(x,y,z)$ and accompanying upper bounds $U$. Similar results for various $\Phi$ in the Lorenz system appear in \cite{Goluskin2017}, along with a more detailed discussion of computational implementation. Here we report more precise computations for $\Phi=z^4$, obtained using the multiple precision SDP solver SDPA-GMP \cite{Yamashita2010, Nakata2010}. Conversion of SOS conditions to SDPs was automated by YALIMP \cite{Lofberg2004, Lofberg2009}, which was interfaced with the solver via mpYALMIP \cite{mpyalmip}.

\setlength{\tabcolsep}{4pt}
\begin{table}
\centering
\caption{\label{tab: bounds} Upper bounds $\ov{z^4}\le U$ in the Lorenz system computed using polynomial $V(x,y,z)$ of various degrees. Underlined digits agree with the value $\ov{z^4}\approx 592827.338$ attained on the shortest periodic orbit.}
\begin{tabular}{cl}
Degree of $V$ & Upper bound $U$ \\
\hline
4 & \quad 635908. \\
6 & \quad\underline{59}5152. \\
8 & \quad\underline{592}935. \\
10 & \quad\underline{592827}.568  \\
12 & \quad\underline{592827.3}44
\end{tabular}
\end{table}

Table \ref{tab: bounds} reports upper bounds computed by solving SDPs that produce optimal $V$ of various polynomial degrees. As the degree of $V$ increases, the bounds approach the value of $\ov{z^4}$ on the shortest periodic orbit to within 7 significant figures. This suggests that $\ov{z^4}$ is maximized on this orbit, and it reflects the sharpness of the bounds asserted by the duality \eqref{eq: strong}. We do not report the lengthy expressions for these $V$; some simpler examples appear in~\cite{Goluskin2017}.

To demonstrate how the volumes $\mathcal{S}_M$ defined in \eqref{eq: S} approximate maximal trajectories, we consider the polynomials $V$ of degrees 6 and 10 that produce the bounds in Table \ref{tab: bounds}. For the maximum in the definition of $\mathcal{S}_M$ we use the corresponding $U$, which bounds it from above. In each case we find that $U$ is within $0.1$ of the true maximum over any ball $B$ containing the \mbox{attractor}.

Figure \ref{fig: 3D}a shows the volume $\mathcal S_{3000}$ for the degree-6 $V$, as well as the orbit that appears to maximize $\ov{z^4}$. The volume captures the rough location and shape of the orbit while omitting much of the strange attractor, but this $V$ is not optimal enough to yield strong quantitative statements. It follows from \eqref{eq: Cheb est} that any trajectory where $\ov{z^4}$ is within $\epsilon$ of the upper bound $U = 595152$ must lie inside $\mathcal S_{3000}$ for a fraction of time no less than $1-\epsilon/3000$. However, there are no trajectories on which this is close to unity; the higher-degree bounds in Table \ref{tab: bounds} preclude any trajectories with $\epsilon<2324$.

The degree-10 $V$ gives a significantly refined picture of maximal and nearly maximal trajectories for $\ov{z^4}$. Figure \ref{fig: 3D}b shows the volume $\mathcal S_{1000}$ defined using this $V$.
It follows from \eqref{eq: Cheb est} that any trajectory where $\ov{z^4}$ comes within $\epsilon$ of $U = 592827.568$ must lie in $\mathcal S_{1000}$ for a fraction of time no less than $1-\epsilon/1000$. There exist trajectories on which this is nearly unity: on the shortest periodic orbit $\ov{z^4}$ is only $\epsilon\approx 0.23$ smaller than $U$. Any trajectory where $\ov{z^4}$ is so large must spend at least 99.97\% of its time in $\mathcal S_{1000}$.

Finding maximal trajectories directly may be intractable in many systems. We propose that the next best option is to compute volumes like those in Fig.\ \ref{fig: 3D}. However, we caution that finding points in a set $\mathcal S_M$ defined by \eqref{eq: S} can itself be difficult, even for polynomials.

\begin{figure} 
\centering
\begin{tikzpicture}
\node at (0,0) {\includegraphics[width=.7\columnwidth]{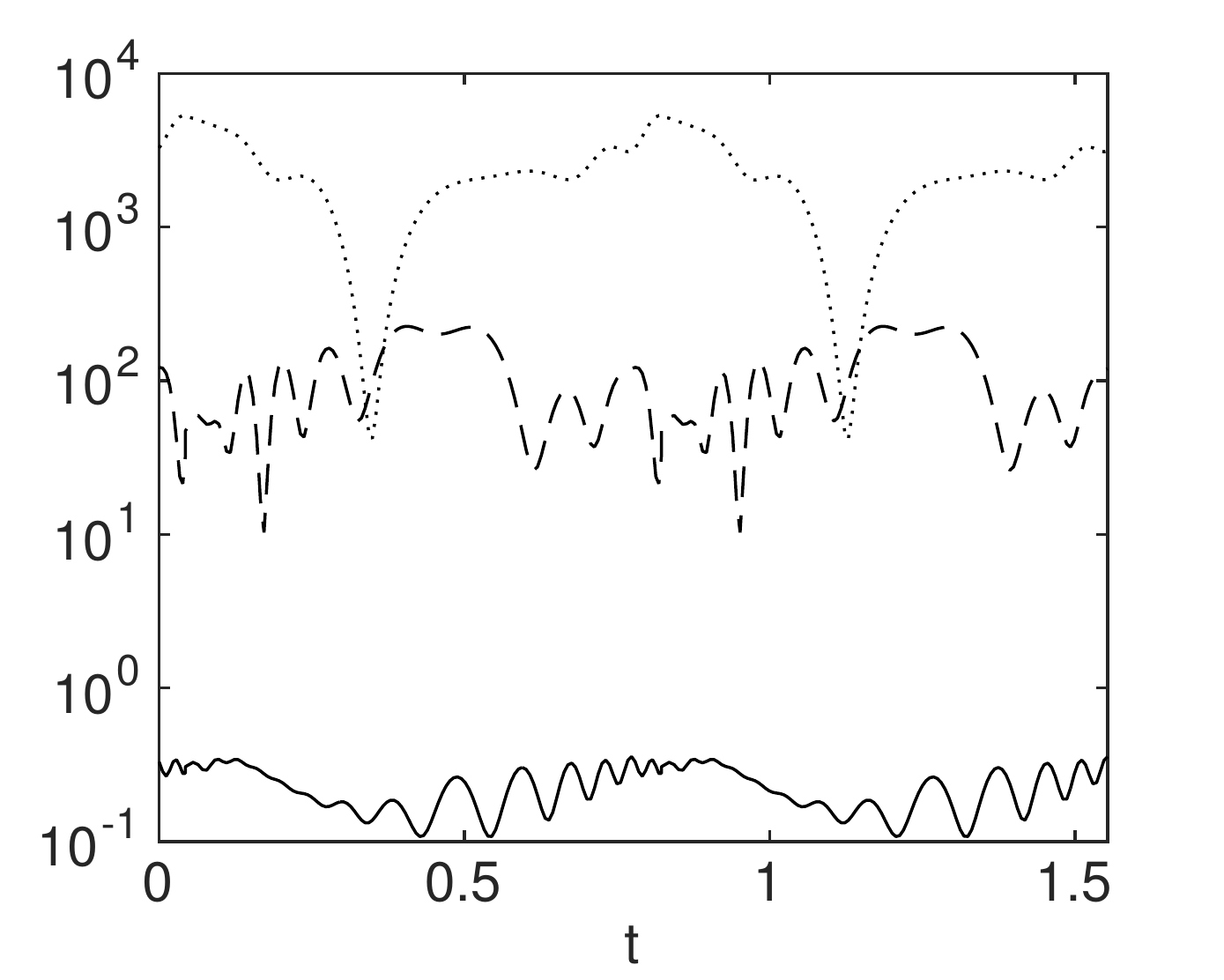}};
\node at (-3.35,2.05) {(a)};
\end{tikzpicture}
\begin{tikzpicture}
\node at (0,0) {\includegraphics[width=.7\columnwidth]{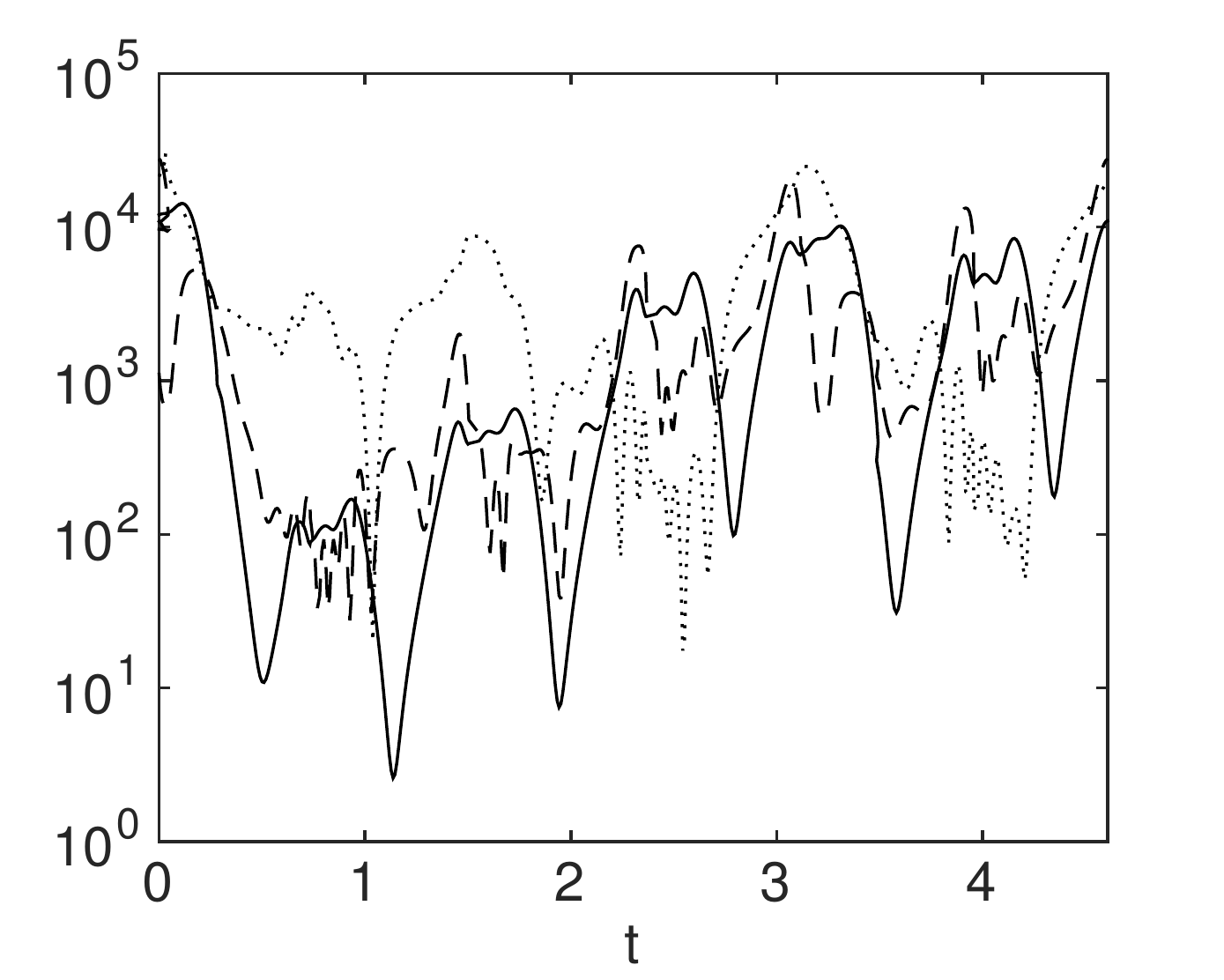}};
\node at (-3.35,2.05) {(b)};
\end{tikzpicture}
\caption{\label{fig: along orbits}The quantity $U - \Phi -\f\cdot\nabla V$ for $\Phi=z^4$ and polynomials $V$ of degrees 6 (\protect\dottedrule), 8 (\protect\dashedrule), and 10 (\solidrule), plotted along (a) the shortest periodic orbit and (b) the periodic orbit with symbol sequence $AABABB$.}
\end{figure}

As the auxiliary functions producing upper bounds on $\ov\Phi^*$ approach optimality, the integrand in \eqref{eq: near optimal ineq} approaches zero almost everywhere on maximal trajectories. This can be seen in Fig.\ \ref{fig: along orbits}a, where the integrand is plotted along the shortest periodic orbit in the Lorenz system for our polynomials $V$ of degrees 6, 8, and 10. Along other orbits where $\ov{z^4}$ is large but not maximal, $V$ is less strongly constrained. As an example, we consider the periodic orbit computed in \cite{Viswanath2003} that winds around the two wings of the Lorenz attractor with symbol sequence $AABABB$. On this orbit $\ov{z^4}$ is smaller than the maximum by approximately 2798. The integral in \eqref{eq: near optimal ineq} remains between 0 and 2798 as $V$ approaches optimality but need not approach 0 on this orbit. In our computations it does not, as seen in Fig.~\ref{fig: along orbits}b.

Although the auxiliary polynomials $V$ yielding the bounds on $\ov{z^4}$ in Table \ref{tab: bounds} approach optimality, they are not exactly optimal. Optimal $V^*$ which are polynomial have been constructed to prove sharp bounds on other averages in the Lorenz system, including $\ov{z}$, $\ov{z^2}$, and $\ov{z^3}$ \cite{Malkus1972, Knobloch1979, Goluskin2017}. These averages are maximized on the two nonzero equilibria; in each case the set $\mathcal S_0$ corresponding to $V^*$ is the line through these equilibria. These $\mathcal S_0$ notably include points not on any maximal trajectory. In contrast, for $\ov{z^4}$ the shortest periodic orbit appears to be maximal. This conjecture could be proved by constructing a $V$ whose $\mathcal S_0$ contains the shortest orbit. If such a $V$ exists, it would necessarily be optimal. However, we expect that this orbit is non-algebraic and that no polynomial $V$ can be optimal.

\smallskip
\noindent
{\bf 5. Proof of duality} \quad
To prove the strong duality \eqref{eq: strong} we require several facts from ergodic theory, which are provable by standard methods as in \cite{walters1982introduction}. (See also \cite[Chap.\ 12]{phelps2001lectures}.) Let $\varphi_t(\x)$ denote the flow map $\x(\cdot)\mapsto\x(\cdot+t)$ for the ODE \eqref{eq: ode}. By assumption, $\varphi_t$ is well-defined on $B$ for all $t\geq0$ and is continuously differentiable there. Let $Pr(B)$ denote the space of Borel probability measures on $B$. A measure $\mu\in Pr(B)$ is \emph{invariant} with respect to $\varphi_t$ if $\mu(\varphi_t^{-1}A) = \mu(A)$ for all Borel sets $A$ and all $t$. Such a measure is \emph{ergodic} if to any invariant Borel set it assigns measure either zero or one. The set of invariant probability measures on $B$ is nonempty, convex, and weak-$*$ compact; its extreme points are ergodic.

Our proof of the duality \eqref{eq: strong} proceeds via a standard minimax template from convex analysis (see, e.g., \cite{ekeland1999convex}). It suffices to establish the following sequence of equalities:
\begin{subequations}
\begin{align}
\max_{\x_0\in B}\,\overline\Phi
\label{eq: relaxed} &= \max_{\substack{\mu\in\Pr(B) \\ 
\mu\text{ is invar.}}}\,\int \Phi\,d\mu \\
\label{eq: sup-inf} &= \sup_{\mu\in Pr(B)}\inf_{V\in C^{1}(B)}\, \int \Phi+\mathbf{f}\cdot\nabla V\,d\mu \\
\label{eq: inf-sup} &= \inf_{V\in C^{1}(B)}\sup_{\mathbf{\mu}\in Pr(B)}\, \int \Phi+\mathbf{f}\cdot\nabla V\,d\mu \\
\label{eq: trivial} &= \inf_{V\in C^{1}(B)}\max_{\x\in B}\,\left\{ \Phi+\mathbf{f}\cdot\nabla V\right\}.
\end{align}
\end{subequations}
In \eqref{eq: relaxed} we reformulate our problem as a maximization over invariant measures, whose analogue for discrete maps is the topic of the field of ergodic optimization \cite{Jenkinson2006}. The remainder of this section is devoted to proving the first three equalities \eqref{eq: relaxed}--\eqref{eq: inf-sup}, along with the fact that the maximum in \eqref{eq: relaxed} is attained. The final equality \eqref{eq: trivial} is evident since, for each $V$, the supremum in \eqref{eq: inf-sup} is attained by a suitable Dirac measure.

We begin by proving \eqref{eq: relaxed}. We claim that the right hand problem appearing there is a concave relaxation of the lefthand problem, and that it attains the same maximum. To see this, note first that for each initial condition $\x_0$ in $B$ there exists an invariant probability measure $\mu$ that attains $\ov\Phi(\x_0) = \int \Phi \,d\mu$. Thus,
\begin{equation}
\sup_{\x_0\in B}\,\ov\Phi(\x_0) \leq \max_{\substack{\mu\in\Pr(B)\\ 
\mu\text{ is invar.}}}\int \Phi\,d\mu.
\label{eq:relaxation}
\end{equation}
The righthand problem in \eqref{eq:relaxation} is a maximization of a continuous linear functional over a compact convex subset of $Pr(B)$, so it achieves its maximum at an extreme point $\mu^*$ \cite[Chap.\ 13]{lax2002functional}, which is an ergodic invariant measure.
By Birkhoff's ergodic theorem \cite{walters1982introduction},  
\begin{equation}
\overline{\Phi}(\mathbf{x}_0) =
\int \Phi \,d\mu^* = 
\max_{\substack{\mu\in\Pr(B)\\ 
\mu\text{ is invar.}}}\int \Phi\,d\mu
\end{equation}
for almost every $\x_0$ in the support of $\mu^*$.
Therefore the inequality in \eqref{eq:relaxation} is in fact an equality, and any such $\x_0$ attains the maximal time average $\ov\Phi^*$. This proves \eqref{eq: relaxed}.

To prove the second equality \eqref{eq: sup-inf} we require the following equivalence of Lagrangian and Eulerian notions of invariance: a Borel probability measure $\mu$ is invariant with respect to $\varphi_t$ by the usual (Lagrangian) definition if and only if the vector-valued measure $\mathbf{f}\mu$ is weakly divergence-free. The latter condition, which we denote by $\text{div}\,\mathbf{f}\mu=0$, means that
\begin{equation}
\int \mathbf{f}\cdot\nabla\psi\,d\mu=0
\label{eq: Eul inv}
\end{equation}
for all smooth and compactly supported $\psi(\x)$.
This is an Eulerian characterization of invariance.

The fact that $\text{div}\,\mathbf{f}\mu=0$ is equivalent to invariance is quickly proved using the flow semigroup identity, which states that $\varphi_{t+s}=\varphi_{t}\circ\varphi_{s}$ for all $t$ and $s$.
It follows that
\begin{equation}
\frac{d}{dt}\int \psi\circ \varphi_t\,d\mu = \int \mathbf{f}\cdot \nabla (\psi\circ \varphi_t)\,d\mu
\label{eq: invariance_equiv}
\end{equation}
for all smooth and compactly supported $\psi$. If $\text{div}\, \mathbf{f}\mu=0$, the righthand side of \eqref{eq: invariance_equiv} vanishes, so $\mu$ is invariant. Conversely, if $\mu$ is invariant then the lefthand side of \eqref{eq: invariance_equiv} vanishes for all $t$, and at $t=0$ we find the statement that $\mathbf{f}\mu$ is weakly divergence-free.

With the Eulerian characterization of invariance in hand, we turn to proving \eqref{eq: sup-inf}.
Depending on $\mu$, there are two possibilities for the minimization over $V$ in \eqref{eq: sup-inf}:
\begin{equation}
\inf_{V\in C^{1}(B)}\,\int \mathbf{f}\cdot\nabla V\,d\mu = \begin{cases}
0 & \text{div}\, \f\mu = 0 \\
-\infty & \text{otherwise}.
\end{cases}
\end{equation}
Only measures for which $\text{div}\, \f\mu=0$ can give values larger than $-\infty$ in \eqref{eq: sup-inf}. As shown above, $\text{div}\, \f\mu=0$ if and only if $\mu$ is invariant. Therefore, since there always exists at least one invariant probability measure, 
\begin{equation}
\sup_{\mu\in Pr(B)}\inf_{V\in C^{1}(B)}\, \int \Phi+\mathbf{f}\cdot\nabla V\,d\mu
=\max_{\substack{\mu\in\Pr(B)\\
\mu\text{ is invar.}
}}
\int\Phi\,d\mu.
\label{eq: sup-inf 3}
\end{equation}
Thus \eqref{eq: sup-inf} is proven. In other words,
\begin{equation}
\mathcal{L}(\mu,V)=\int\Phi+\f\cdot\nabla V\,d\mu
\label{eq: Lagrangian}
\end{equation}
is a Lagrangian for the constrained maximization appearing on the righthand side of \eqref{eq: sup-inf 3}.

Finally, we prove the equality \eqref{eq: inf-sup}. In terms of the Lagrangian $\mathcal{L}$, we must show that
\begin{equation}
\sup_{\mu\in\Pr(B)}\inf_{V\in C^{1}(B)}\, \mathcal{L} = \inf_{V\in C^{1}(B)}\sup_{\mu\in\Pr(B)}\, \mathcal{L}.
\label{eq: infsup_supinf}
\end{equation}
The fact that the order of $\inf$ and $\sup$ can be reversed without introducing a so-called duality gap is not trivial; it is at the heart of our proof of the strong duality \eqref{eq: strong}.
This reversal relies on properties of the Lagrangian $\mathcal L$ and the spaces $Pr(B)$ and $C^1(B)$.

The desired equality \eqref{eq: infsup_supinf} can be proved using any of several abstract minimax theorems from convex analysis. Here we apply a fairly general infinite-dimensional version due to Sion \cite{sion1958general}. We follow the notation of its statement in the introduction of \cite{komiya1988elementary}, which contains an elementary proof. Let $X=\Pr(B)$ in the weak-$*$ topology. It is a compact convex subset of a linear topological space. Let $Y=C^1(B)$ in the $C^1$-norm topology, which is itself a linear topological space. Take $f = - \mathcal{L}$ and observe that $f(x,\cdot)$ is upper semicontinuous and quasi-concave on $Y$ for each $x\in X$, and that $f(\cdot, y)$ is lower semicontinuous and quasi-convex on $X$ for each $y\in Y$. Then \eqref{eq: infsup_supinf} follows from a direct application of Sion's minimax theorem \cite{komiya1988elementary}, so \eqref{eq: inf-sup} is proven. 

This completes the proof of the equalities \eqref{eq: relaxed}--\eqref{eq: trivial} and so too the proof of the strong duality~\eqref{eq: strong}.

\smallskip
\noindent
{\bf 6. Conclusions} \quad
This Letter establishes that the auxiliary function method for proving \emph{a priori} bounds on long-time averages in dynamical systems yields arbitrarily sharp bounds, so long as the dynamics arise from ODEs. The proof elucidates the role that auxiliary functions play in the search for optimal bounds: they are Lagrange multipliers enforcing the constraint of invariance for probability measures on phase space. We also have demonstrated that certain sets constructed from nearly optimal auxiliary functions can be used to locate all optimal and nearly optimal trajectories. How close an auxiliary function is to optimality determines the fraction of time nearly optimal trajectories are guaranteed to spend in these sets. We expect much can be learned about the shape of optimal trajectories and their invariant measures by the auxiliary function approach.

Many of these observations extend to infinite-dimensional dynamics that are governed by nonlinear partial differential equations (PDEs) of the form 
\begin{align}
\frac{d}{dt} \mathbf{u} = \mathbf{f} \{\mathbf{u} \}.
\label{eq: PDE}
\end{align}
Auxiliary functionals $V\{\mathbf{u}\}$ defined on a suitable function space yield \emph{a priori} bounds on long-time averages just as in the finite-dimensional case. The ``background method'' used to bound mean quantities in fluid dynamics and other systems \cite{Doering1992} is an example of using quadratic $V$ \cite{Chernyshenko2017}. Whether or not nearly optimal functionals are always guaranteed to exist, and if they are ever quadratic for systems of interest, remains unclear. This emphasizes the need for a rigorous proof of duality between auxiliary functionals and extremal trajectories for general PDEs. 

\smallskip
\noindent
{\bf Acknowledgements} \quad
We thank Lora Billings, Rich Kerswell, Edward Ott, Ralf Spatzier, Divakar Viswanath, and Lai-Sang Young for helpful discussions and encouragement.
This work was supported by NSF Award DMS-1515161, Van Loo Postdoctoral Fellowships (IT, DG), and a Guggenheim Foundation Fellowship (CRD).

\bibliographystyle{apsrev4-1}
\bibliography{bounding_duality_v32.bbl}

\end{document}